\documentclass[10pt, twocolumn, final]{IEEEtran}  
\IEEEoverridecommandlockouts  

\usepackage{booktabs} 

\usepackage{color,verbatim}
\definecolor{red}{rgb}{1,0.2,0.2}
\definecolor{green}{rgb}{0.2,1,0.5}
\definecolor{blue}{rgb}{0,0,1}
\definecolor{lightblue}{rgb}{0.3,0.5,1}

\usepackage[cmex10]{amsmath}
\usepackage{amsfonts}
\usepackage{amssymb}
\usepackage{verbatim} 
\usepackage{algorithm}
\usepackage{algorithmic}
\usepackage{graphics,psfrag,epsfig}
\usepackage{epstopdf} 
\usepackage{hyperref}

\newtheorem{theorem}{Theorem}

\newtheorem{corollary}{Corollary}

\newtheorem{rem}{Remark}

\newcommand{\nnb}{\nonumber}
\newcommand{\1}{\mathbf{1}}
\newcommand{\0}{\mathbf{0}}
\newcommand{\ab}{\mathbf{a}}

\newcommand{\eb}{\mathbf{e}}

\newcommand{\T}{^{\mathsf{T}}}

\newcommand{\yb}{\mathbf{y}}

	\title{\LARGE Optimal Cache Allocation for Named Data Caching under Network-Wide Capacity Constraint}
\author{Van Sy Mai, Stratis Ioannidis, Davide Pesavento, and Lotfi Benmohamed%
	\thanks{
		V.~S.~Mai, D.~Pesavento and L.~Benmohamed are with the Information Technology Laboratory, NIST, USA. Emails: {\tt \{vansy.mai,
			davide.pesavento,lotfi.benmohamed\}@nist.gov}. \mbox{S.~Ioannidis} is with the Northeastern University and is supported by NSF grants NeTS-1718355 and CCF-1750539. Email: {\tt ioannidis@ece.neu.edu}.
		
	 Mention of commercial products does not imply NIST's endorsement.}%
}
\begin{document}	
	\maketitle	
	
	\begin{abstract}
	Network cache allocation and management are important aspects of the design of an Information-Centric Network (ICN), such as one based on Named Data Networking (NDN). We address the problem of optimal cache size allocation and content placement in an ICN in order to maximize the caching gain resulting from routing cost savings. While prior art assumes a given cache size at each network node and focuses on content placement, we study the problem when a global, network-wide cache storage budget is given and we solve for the optimal per-node cache allocation. This problem arises in cloud-based network settings where each network node is virtualized and housed within a cloud data center node with associated dynamic storage resources acquired from the cloud node as needed. With the offline centralized version of the optimal cache allocation problem being NP-hard, we develop a distributed adaptive algorithm that provides an approximate solution within a constant factor from the optimal. Performance evaluation of the algorithm is carried out through extensive simulations involving a variety of network topologies, establishing experimentally that our proposal significantly outperforms existing cache allocation algorithms.
	\end{abstract}

	\begin{IEEEkeywords} Caching, distributed optimization, Information-Centric Networking
	\end{IEEEkeywords}

	\section{Introduction}
	
	Traditional networking is being transformed into a more agile one with significant flexibility in how network services get deployed. Networking hardware is being virtualized by a cloud infrastructure (using hypervisor/virtualization software) with network forwarding done in software within virtual machines (VMs). Software-defined networking (SDN) and network functions virtualization (NFV) are among the enablers of this virtualization of networking. The telecom operator's business is evolving and traditional network operators are increasingly becoming ``Telecom Cloud Operators''. These operators are deploying their own cloud infrastructure with dedicated data centers to meet multiple objectives: (a) for the deployment of their own telecom network in support of their core telecom business, (b) in support of their own IT needs, and (c) to get into the cloud market currently served by cloud operators. While the cloud-based network (CBN) resulting from objective (a) above is built using the operator's private cloud, large companies that want to build their private wide-area enterprise network as a CBN can do so using resources from public cloud operators.
	
	When networks are virtualized (whether private of public CBNs), they become more flexible and dynamic in many aspects including in their caching capability. When these networks are deployed using ICN technology, they will be able to implement a dynamic cache feature. The network operator will now be able to size its ICN's cache dynamically to maintain good performance as the produced content changes and associated popularity evolves with the user demand profile. Given a network-wide cache budget $M$ that the operator is willing to invest, the problem to address in this case is how should this budget be dynamically allocated among the network nodes to maximize performance under varying network conditions. Shifting cache capacity among network nodes over time will be an easy task within a CBN: the per-node cache storage capacity can be increased or decreased as needed (by acquiring or relinquishing storage blocks from/to the storage pool at the cloud site where the node is homed) while staying within the preset network-wide limit $M$. While the problem of assigning items to caches under given fixed cache sizes has already been studied, this cache capacity design problem has not, and it is the focus of this paper.
	
	In this paper we address the modeling, analysis, and implementation of caching in cloud-based information centric networks. In these networks a subset of nodes act as the designated sources for content (data producers) while any node can be a data consumer that generates requests for data items, which get forwarded toward the designated producers. These requests may not reach the ultimate producer as ICN forwarding ends when reaching a node along the path that has cached the requested item in its Content Store (CS). When such a cache hit takes place, the requested item is served from the CS and sent back to the requesting node along the reverse path.
	
	Literature on ICN caching is extensive \cite{carofiglio2015lac,thomas2015object,nguyen2015congestion,chai2012cache,dehghan2016utility,ming2014age,badov2014congestion,ioannidis2016adaptive}. With an ICN being a network of caches where each network node is equipped with a content store, designing a good caching solution involves the aspects of determining the size of each CS, deciding which data objects should be cached (placement strategy), and which ones should be evicted when needed (replacement strategy). An efficient caching solution brings many benefits as it (a) reduces the data producer load since consumer's requests would rarely be satisfied by the producer but rather by cashes, (b) significantly reduces the amount of network traffic and avoids bottlenecks caused by publishing data at a limited set of locations, and (c) offers users a faster content retrieval for an enhanced user experience. In other words, the investment in caching is expected to be of benefit to users, network operators, as well as content providers when it enables performance similar to content distribution networks (CDNs) by dynamically storing content in regions of high demand.
	
	Our goal is to achieve an optimal caching solution that maximizes the caching gain by minimizing the aggregate routing costs due to content retrieval across the network. The network load made up of each user demand, which is determined by the rate of requests and the paths they follow, is typically dynamic and not known in advance. As a result it is desirable to have adaptive caching solutions that can achieve optimal placement of data items in network caches without prior knowledge of the demands and be adaptive to any potential demand changes. In addition to being adaptive, caching needs to be distributed as well, since centralized solutions are not expected to be feasible when multiple administrative domains are involved. The network is expected to be more scalable when implementing distributed algorithms with caching decisions that rely only on locally available information.
	
	Path replication, also known as Leave Copy Everywhere (LCE), is a popular caching strategy that is dynamic and distributed, and is often discussed in the literature \cite{che2002hierarchical,jacobson2009networking,lv2002search,rossi2011caching}. When a data item is forwarded on the reverse path towards the consumer that requested it, it is cached at each intermediate node along the path. When a node's cache is full, a replacement takes place by evicting an already cached item using policies such as LRU, LFU, or FIFO. Despite its popularity, LCE has no performance guarantees and can be shown to be arbitrarily suboptimal~\cite{ioannidis2016adaptive}.
	
	In this paper, we discuss our design of a distributed and adaptive caching solution with provable performance guarantees. Our main contributions are the following:
	\begin{itemize}
		\item While previous work with a similar problem formulation assumes that cache sizes are given and only deals with object placement in such fixed size caches, we address a more general problem where no assumption of fixed cache sizes is made but rather uses a global network-wide cache budget constraint, and design the optimal per-node cache capacity.
		\item We make use of a game theory framework to design a distributed algorithm for this problem by combining a distributed gradient estimation with the distributed constraint satisfaction methodology.
		
		\item We show that our game-based algorithm can provide suboptimal solutions within a factor $(1 - 1/e - \epsilon)$ of optimum for any given small $\epsilon > 0$ and without prior knowledge of the network demand. 
		
		\item We present results from extensive simulations over a number of network topologies that show how our algorithm outperforms those based on fixed size caches.
	\end{itemize}
	
	The remainder of this paper is organized as follows. In Section~\ref{sec:related} we review related work. We introduce the system model and formally state the problem in Section~\ref{sec:problem}. Centralized solutions to the problem are discussed in Section~\ref{sec:central}. Our main results on distributed algorithms are discussed in Section~\ref{sec:distr} along with a discussion on implementation issues. Numerical results are presented in Section~\ref{sec:sim} and followed by conclusions in Section~\ref{sec:conclusion}.

	\section{Related Work}
	\label{sec:related}
	
	The offline problem we study amounts to maximizing a submodular function subject to matroid constraints. Such problems are ubiquitous and appear in many domains (see Krause and Golovin \cite{krause2012} for a detailed overview). Though NP-hard, there exist known approximation algorithms:  Nemhauser et al.~\cite{greedy2} show that the greedy algorithm produces a solution  within 1/2 of the optimal.   Vondr\'ak \cite{vondrak2008optimal} and Calinescu et al.~\cite{calinescu2007maximizing,calinescu2011} show that the so-called continuous-greedy algorithm produces a solution within $(1-1/e)$ of the optimal in polynomial time, which cannot be further improved \cite{nemhauser1978best}. The latter requires sampling the so-called multi-linear relaxation of the objective.
	
	Specifically in the context of caching gain maximization, Shanmugam et al.~\cite{shanmugam2013femtocaching} and Ioannidis and Yeh \cite{ioannidis2016adaptive,ioannidis2018adaptive} consider a more restricted version of our problem, in which (a) cache sizes are given, and (b) only object placements are optimized. Shanmugam et al.~study this under a restricted topology, assuming homogeneous (i.e., equal-size) caches; Ioannidis and Yeh study the problem under arbitrary topology and cache sizes. The authors show in each of these settings, respectively, that the concave relaxation technique of Ageev and Sviridenko \cite{ageev2004pipage} also attains the $1-1/e$ approximation ratio; this algorithm is preferable to the continuous-greedy algorithm, as it eschews sampling. Ioannidis and Yeh further propose a distributed adaptive algorithm with the same approximation ratio, based on projected subgradient ascent of the relaxation function $L$ (see Section~\ref{secCVX_Relaxation}). A similar approach can be used to jointly optimize both caching and routing decisions \cite{ioannidis2017icn,ioannidis2018jointly}.
	
	In contrast to \cite{ioannidis2016adaptive,ioannidis2018adaptive}, the setting we study does not assume fixed cache sizes. In optimization terms, our problem includes an additional global constraint, introduced through the global budget $M$. Feasible solutions still define a matroid under this additional constraint; moreover, we show that the concave relaxation technique of Ageev and Sviridenko \cite{ageev2004pipage} also applies to this offline setting. However, the global coupling through this constraint makes the projected gradient ascent method of Ioannidis and Yeh inadequate as a distributed adaptive algorithm. 
	Our problem also resembles network resource allocation or utility maximization problems, where various decomposition techniques admit  distributed implementations \cite{palomar2006tutorial}. However, our (relaxed) global cost function is coupled in a such way that renders decomposition and decoupling approaches communication-expensive
	and complicated, let alone distributed adaptive implementation. This calls for a different approach. 
	Specifically, we propose to employ the game theory-based framework in  \cite{marden2012state,li2013designing,li2014decoupling,marden2015game} for designing a distributed algorithm, where the global cost is embedded in the potential of a game. In contrast to \cite{marden2012state,li2013designing,li2014decoupling,marden2015game}, however, we do not assume separability of the potential function, nor do we employ any decoupling technique; this is achieved by making use of the distributed gradient estimation scheme in \cite{ioannidis2016adaptive}. 
	
	Finally, our work is also related to the problem of  virtual machine (VM) allocation  in  cloud computing \cite{li2011virtual,guenter2011managing,van2010cost,batista2007set}--see also \cite{jiang2012joint}, that jointly optimizes placement and routing in this context. Heterogeneity of host resources and VM requirements leads to multiple knapsack-like constraints (one for each resource) per host. Our storage constraints are simpler; as a result, in contrast to \cite{li2011virtual,guenter2011managing,van2010cost,batista2007set,jiang2012joint},  we can provide distributed algorithms with provable approximation guarantees.

	\section{Preliminaries and Problem Formulation}
	\label{sec:problem}
	
	\subsection{Notational Conventions.}
	In what follows, we denote by $\mathbb{R}$,   $\mathbb{N}$ the sets of real and natural numbers, respectively. If $\mathcal{A}$ is a finite set, $|\mathcal{A}|$ denotes its cardinality. Let $[\cdot]_{\mathcal{X}}$ denote the projection operator onto the set $\mathcal{X}$; $[\cdot]_+$ is the same as $[\cdot]_{\mathbb{R}_+}$.
	
	\subsection{System Model}
	Consider a \emph{connected} network $\mathcal{G} = (\mathcal{V},\mathcal{E})$  where $\mathcal{V}$ is the set of  nodes and $\mathcal{E} \subset \mathcal{V} \times \mathcal{V}$ is the set of links. Nodes are equipped with caches (content stores), whose capacity can be adjusted as part of an optimized design. As discussed above, the nodal cache size can be adjusted as needed by acquiring or relinquishing units of storage at the local cloud node (data center) part of the operator's deployed cloud. 
	The local cache is used to store content items from a catalog made up of a set  $\mathcal{C}$, and  subsequently serve requests for these items from the cache.
	We denote by $M$ the total cache capacity that the network operator is willing to deploy network-wide, it reflects a limit on the operator's budget invested in network cache storage.
	
	We denote by $\bar{c}_v \in \mathbb{N}$ the maximum cache capacity that can be allocated at node $v$; this restriction would typically be due to limits on the available physical storage at the cloud node where the network node $v$ is homed. Note that it is more likely that the network limit $M$ and the available physical storage in the cloud are such that the $cv$ likely will likely not be reached. However, our model and the following analysis and design that follows is capable of handling both nodal-level and network-level capacity limits.
	We denote by  $x_{v,i}\in \{0,1\}$ for $v\in \mathcal{V}$, $i\in \mathcal{C}$ the variable indicating if node $v$ stores item $i$. 
	To store these items, the capacity at node $v$ is thus $\sum_{i\in\mathcal{C}}x_{vi}$, which must be less than $\bar{c}_v$. Moreover, given the budget constraint, the total capacity $\sum_{v\in \mathcal{V}}\sum_{i\in \mathcal{C}}x_{vi}$ must be less than $M$.

	We assume that, for each item $i$, there exists a set $\mathcal{S}_i \subset \mathcal{V}$ of nodes that serve as designated servers for that item (data producers): these nodes always store $i$, i.e., $x_{vi} = 1, \forall v\in \mathcal{S}_i$.
	Requests arrive in the network and traverse predetermined paths towards the designated servers of each item. Formally, a request for item $i\in \mathcal{C}$ through path $p = \{ p_1,\ldots,p_K \}\subset \mathcal{V}$ is denoted by pair $(i,p)$.  We denote by $\mathcal{R}$ the set of all such requests.
	We assume that requests $(i,p)\in \mathcal{R}$ are \emph{well-routed}, i.e., follow paths with no loops that terminate at designated servers in $\mathcal{S}_i$.
	Moreover, requests for each element in $\mathcal{R}$ arrive according to independent Poisson processes with arrival rates $\lambda_{(i,p)}>0$; note that such assumption is standard for modeling request arrivals (e.g., \cite{carofiglio2015lac,nguyen2015congestion,chai2012cache,dehghan2016utility,ming2014age,badov2014congestion,ioannidis2016adaptive,che2002hierarchical}).
	
	A request $(i,p)$ is routed following path $p$ until it reaches a cache that has item $i$.
	Then, a response message carrying item $i$ is generated and sent over $p$ in the reverse direction back to the first node in $p$.
	We assume that the cost of routing an item over a link $(i,j) \in \mathcal{E}$ is $w_{ij} \in \mathbb{R}_+$, while the cost of forwarding requests is negligible.
	
	The goal of the network designer is to select (a) the cache capacity at each node, as well as (b) which items to store at each cache.
	The purpose of our design is to jointly allocate storage resources and item placement at each node in order to minimize routing costs.
	To be adaptive, this allocation should occur dynamically, without assuming prior knowledge of user demand for items. We formalize this optimization problem below.
	
	\subsection{Problem Statement}
	
	Recall that the network designer acquires storage at each network node from the local cloud node subject to a prescribed budget $M$. We seek a joint item placement and cache capacity allocation that minimizes the aggregate expected cost. In particular, let $C_0$ denote the expected cost when there are no items cached except for designated servers, i.e., 
	\begin{equation}\label{eqC_0}
	C_0  =  \sum_{(i,p)\in \mathcal{R}} \lambda_{(i,p)} \sum_{k=1}^{|p|-1} w_{p_{k+1} p_k}.
	\end{equation} 
	In the presence of cached contents, the cost of serving a request $(i,p)\in \mathcal{R}$ is
	\begin{equation}\label{eqC_ip}
	C_{(i,p)}(X)  =  \sum_{k=1}^{|p|-1} w_{p_{k+1} p_k}\prod_{l=1}^k (1-x_{p_{l}i}).
	\end{equation} 
	Thus,  the expected caching gain corresponding to an allocation $X = \{x_{vi}\}_{v\in \mathcal{V}, i\in \mathcal{C}}$ is given by $F(X) := C_0 - \sum_{(i,p)\in \mathcal{R}} \lambda_{(i,p)} C_{(i,p)}(X)$, i.e., 
	\begin{align}
	F(X)  = \sum_{(i,p)\in \mathcal{R}} \lambda_{(i,p)} \sum_{k=1}^{|p|-1} w_{p_{k+1} p_k} \Big( 1-\prod_{l=1}^k (1-x_{p_{l}i}) \Big). \label{eqF_X}
	\end{align} 
	Formally, we seek to develop a distributed adaptive algorithm for the following problem:

	{\it   Given a cache budget $M$ for the whole network, design an allocation $X$ so as to maximize the expected caching gain:
		\begin{align}
		\emph{(MaxCG)} \quad \max_{X} \quad &F(X) & & \nnb\\
		\mathrm{s.t.} \quad  
		&x_{vi} \in \{0,1\}, &\forall& v\in \mathcal{V}, \forall  i \in \mathcal{C}  \label{BinaryConstraint}\\
		&x_{vi} = 1, &\forall& v\in \mathcal{S}_i, \forall i\in \mathcal{C} \label{DesignatedSources}\\
		& \sum_{i\in \mathcal{C}} x_{vi} \le  \bar{c}_v,  &\forall& v\in \mathcal{V} \label{LocalCap} \\
		& \sum_{v\in \mathcal{V}}\sum_{i\in \mathcal{C}} x_{vi} \le M, & &\label{GlobalCap}
		\end{align}
		where $\{\bar{c}_v\}_{v\in \mathcal{V}}$ are nodal maximum allowable capacities.}
	
	\noindent
	For convenience, let $\mathcal{D}_1$ denote the feasible set of (MaxCG), i.e., 
	\begin{equation}
	\mathcal{D}_1 = \{ X\in \mathbb{R}^{|\mathcal{V}|\times|\mathcal{C}|} ~|~\eqref{BinaryConstraint}-\eqref{GlobalCap} \text{~hold} \}.
	\end{equation}
	
	\noindent
	We seek both \emph{centralized}, \emph{offline} algorithms for this problem, as well as \emph{distributed}, \emph{adaptive} algorithms. In the former case, the designer has full knowledge of the demand (request arrival rates) and a full view of system parameters, and solves the problem offline. In the latter case, the caches themselves should adapt and update their capacities to solve the underlying optimization problem in a distributed fashion, by exchanging appropriate messages.

	\begin{rem}    
		A special case is when $\bar{c}_v \ge |\mathcal{C}|$ for all $v\in \mathcal{V}$, that is, local constraints \eqref{LocalCap} are redundant. Another special case is when $M \ge \sum_{v\in \mathcal{V}} \bar{c}_v$, then constraint \eqref{GlobalCap} is redundant and the problem reduces to that in \cite{ioannidis2016adaptive}, which only considers decoupled constraint \eqref{LocalCap}. This paper will focus on the case $M < \sum_{v\in \mathcal{V}} \bar{c}_v$. 
		Thus, our problem has both coupled constraints and coupled objective, and the algorithm in \cite{ioannidis2016adaptive} is no longer applicable. In particular, constraint \eqref{GlobalCap} induces a \emph{global} constraint, coupling decisions \emph{throughout} the network. Maintaining this constraint thoughout the network in a distributed fashion introduces a challenge not present in \cite{ioannidis2016adaptive}, as we discuss below.
	\end{rem}

	\section{Centralized Approaches}
	\label{sec:central}
	Since (MaxCG) is NP-hard \cite{ioannidis2016adaptive}, we seek polynomial-time approximation algorithms for its solution.
	
	\subsection{Greedy Algorithms}
	Since $F$ in \eqref{eqF_X} is a nonnegative, monotone and submodular function and the feasible set of (MaxCG) corresponds to a matroid constraint, the standard greedy algorithm generally yields $\frac{1}{2}$ approximation guarantee\footnote{Without the local constraints \eqref{LocalCap} (or equivalently, $\bar{c}_v > |\mathcal{C}|, \forall v\in \mathcal{V}$), the standard greedy algorithm results in a suboptimal solution within a factor $1-e^{-1}$ of the optimal value.}. A ratio of $1-e^{-1}$ can be achieved by the \textit{continuous greedy algorithm} \cite{calinescu2011maximizing}, or by combining the standard greedy algorithm with a \textit{non-oblivious local search} \cite{filmus2012tight}. We adopt a different approach that will provide some insight on how to construct a distributed, adaptive algorithm.

	\subsection{Convex Relaxation}\label{secCVX_Relaxation}
	We first employ the approximation approach of \cite{ioannidis2016adaptive}: this convexifies both the constraint set and the cost function. In particular, the approximation algorithm proceeds as follows:
	\begin{itemize}
		\item  First, use a continuous relaxation for Boolean variables, i.e., 
		\begin{equation}
		{\max_{Y}} \quad 
		\{  F(Y) ~|~ Y \in {\mathcal{D}}_2   \},
		\end{equation}
		where ${\mathcal{D}}_2 = \mathrm{conv}(\mathcal{D}_1)$, i.e., 
		\begin{equation}
		\mathcal{D}_2 = \{ X \in [0,1]^{|\mathcal{V}|\times|\mathcal{C}|} ~|~ \eqref{DesignatedSources}-\eqref{GlobalCap} \text{~hold} \}.
		\end{equation}
		
		\item Then, approximate the non-concave $F$ by 
		\begin{equation} 
		\!\! L(Y) \!=\!\!\! \sum_{(i,p)\in \mathcal{R}} \lambda_{(i,p)} \!\!\sum_{k=1}^{|p|-1}\!\! w_{p_{k+1} p_k} \min\{1, \sum_{l=1}^k y_{p_{l}i} \}\label{eqL}
		\end{equation}
		which is concave and satisfies $L(X) = f(X), \forall X\in \mathcal{D}_1$. Moreover, 
		\begin{equation} \label{eqBound_F} 
		(1-e^{-1})L(Y) \le F(Y) \le L(Y), \quad \forall Y \in {\mathcal{D}}_2. 
		\end{equation} 
		Thus, the resulting problem
		\begin{equation}
		{\max_{Y}} \quad 
		\{  L(Y) ~|~ Y \in {\mathcal{D}}_2 \} \label{ProbCVX}
		\end{equation}
		is convex; in fact, it can be converted into a linear program. Thus an optimal solution $Y^*$ can be computed in strongly polynomial time. The optimal value of \eqref{ProbCVX} is denoted by $L^*$. 
		
		\item Finally, apply the \emph{pipage rounding} technique of Ageev and Sviridenko \cite{ageev2004pipage} to $Y^*$, yielding a suboptimal solution, denoted by $ [Y^*]^{pp}_{\mathcal{D}_1} $, to (MaxCG) that has $1-e^{-1}$ approximation guarantee; see \cite{ioannidis2016adaptive} for details. 	
		Moreover, since $ L( [Y^*]^{pp}_{\mathcal{D}_1} ) = F([Y^*]^{pp}_{\mathcal{D}_1}) \le F({X}^*) = L({X}^*)\le L(Y^*)$, where $X^*$ is an optimal solution to (MaxCG), we have 
		$$ 
		\frac{ F( [Y^*]^{pp}_{\mathcal{D}_1} )}{F({X}^*)} \ge \frac{L( [Y^*]^{pp}_{\mathcal{D}_1} )}{L(Y^*)}.
		$$
		Thus, we can use the RHS (computed in polynomial time) as another approximation ratio of $[Y^*]^{pp}_{\mathcal{D}_1}$. In practice, this is often better than the theoretical ratio $1-e^{-1}$. 
	\end{itemize}

	\section{Distributed Algorithm based on Potential Game}
	\label{sec:distr}
	This section develops an algorithm for dealing with (MaxCG) based on the convex relaxation approach outlined in the previous section. First, we introduce another simple approximation to the nondifferentiable function $L$ in \eqref{ProbCVX}. Then we show that the game theory framework can be applied to the resulting problem. This game-theoretic approach allows us to adapt both cache capacities as well as content allocations in a distributed fashion. 
	
	\subsection{Continuously differentiable approximation}
	The concave relaxation $L$, given by \eqref{eqL}, is not differentiable. Consider:
	\begin{equation}
	\tilde{L}(Y) := \sum_{(i,p)\in \mathcal{R}} \lambda_{(i,p)} \sum_{k=1}^{|p|-1} w_{p_{k+1} p_k} \mathrm{sat}_{\alpha} \big(\sum_{l=1}^k y_{p_{l}i} \big), \label{eqL3}
	\end{equation} 
	where $\alpha \in (0,1)$ is a small number and 
	\begin{align}
	\mathrm{sat}_{\alpha}(x) := \begin{cases}
	1 & \text{if } x\ge 1 + \frac{\alpha}{2}\\
	1- \frac{(1+ \frac{\alpha}{2}-x)^2}{2\alpha} & \text{if } 1-\frac{\alpha}{2} \le x <1 + \frac{\alpha}{2}\\
	x & \text{if } 0\le x< 1-\frac{\alpha}{2}
	\end{cases}
	\end{align}
	is a lower bound of the function $\min\{1, x\}$ on $\mathbb{R}_+$. 
	Thus, $\tilde{L}$ is a concave lower bound of $L$ and  $\lim_{\alpha \to 0^+} \tilde{L} = L$. 
	Indeed, 
	\begin{equation}\label{eqL3Bound}
	\tilde{L}(Y)  \le L(Y) \le \tilde{L}(Y) + \frac{\alpha}{8}C_0, \quad\forall Y \in \mathcal{D}_2.
	\end{equation} 
	Thus, for sufficiently small $\alpha$, the following problem is a good surrogate for \eqref{ProbCVX}:
	\begin{equation}
	\max_Y \quad\{ \tilde{L}(Y)~|~Y \in {\mathcal{D}}_2\}. \label{ProbCVX_Approx}
	\end{equation}
	
	As a side note, this approximation is useful for the problem considered in \cite{ioannidis2016adaptive}, i.e., (MaxCG) without  the global constraint \eqref{GlobalCap}. The authors employed a distributed subgradient algorithm with a diminishing step size and gain-smoothening to deal with the non-differentiability $L$. Here, $\tilde{L}$ is differentiable with a Lipschitz continuous and bounded gradient. As a result, smooth (asynchronous) optimization algorithms such as distributed projected gradient with a constant step size can be used. This may also provide some insights into why the Greedy Path Replication (asynchronous with constant step size) in  \cite{ioannidis2016adaptive} has a good performance although rigorous analysis was absent. 
	
	Note also that one can consider other alternative approximations; e.g., using $\tanh(x)$ instead of $\mathrm{sat}(x)$ yields another lower bound of $L$ that is strictly increasing, smooth, and strongly concave on $X$. This paper will focus on using $\tilde{L}$ given above, the smoothness and Lipschitz property of which are also valuable for the framework described next. 
	
	\subsection{Potential game design}
	First, we restate \eqref{ProbCVX_Approx} as follows: 
	\begin{align}\label{ProbCVX_Game}
	\displaystyle{\max_{\{\yb_v \in \Omega_v\} }} \quad 
	& \tilde{L}( Y)  \\
	\mathrm{s.t.} \quad  
	& \sum_{v \in \mathcal{V}} (\yb_v\T \1 - c^0_v) \le 0, \label{eqConstraintTotalCacheSize}
	\end{align}
	where $\yb_v\T$ is the $v$-th row of $Y$, $\1$ is a column vector of all ones, $c^0_v$ are constants such that $\sum_{v\in \mathcal{V}} c^0_v  = M$, and 
	\begin{equation*}
	\Omega_v :=  \{ \yb \in [0,1]^{|\mathcal{C}|} ~|~ \sum_{i\in \mathcal{C}} y_{i} \le \bar{c}_v, \quad y_i = 1 \text{ if } v\in \mathcal{S}_i  \}.
	\end{equation*}
	Here, we assume that each node $v$ knows the value $c^0_v$; for example, $c^0_v = M/|\mathcal{V}|$, i.e., the average cache size for the network nodes. 
	
	In the following, we employ the game theory framework in  \cite{marden2012state,li2013designing,li2014decoupling,marden2015game} to design a distributed algorithm for this problem.
	In particular, we will design the state based potential game between the nodes so that they will converge to a pure Nash equilibrium that can be made arbitrarily close to an optimal solution of \eqref{ProbCVX_Game}. The crucial differences between our design and that in \cite{li2013designing,li2014decoupling,marden2015game} are the nodal cost functions and the implementation of the learning algorithm. In particular, in contrast to  \cite{li2013designing,li2014decoupling,marden2015game}, we \emph{do not assume that cost functions are separable across nodes} (indeed, the terms of the objective  \eqref{ProbCVX_Approx} are coupled), \emph{nor do we employ any decomposition technique} (see also Remark~\ref{rem_decomposibility} below).

	\subsubsection{Game model}
		
	We begin by presenting a game played by node caches; the evolution of the game via appropriate dynamics, described below, eventually leads to a solution of \eqref{ProbCVX_Game} in a distributed fashion. 
	
	\begin{enumerate}
		\item[1.]	\textit{State space: } Let $Z = (Y,\eb)$ denote the state of the game, where $\eb = \{e_v\}_{v\in \mathcal{V}}$ and $e_v$ is an error term of node $v$ representing an estimation of $(\yb_{v}\T\1 - c^0_v)$. 
		
		\item[2.]  \textit{Actions: } Each node $v$ has a state-dependent action set $\mathcal{A}_v(Z)$, where an action $\ab_v$ is a tuple $\ab_v = (\hat{\yb}_v, \{\hat{e}_{v\to u}\}_{u \in \mathcal{N}_v})$, where $\hat{e}_{v\to u}$ represents the estimate error that node $v$ sends to a direct neighbor $u$, and $\mathcal{N}_v$ denotes the set of node $v$'s neighbors. 
		
		\item[3.] \textit{State dynamics: } For any state $Z = (Y,\eb)$ and action $\{\ab_v\}$, the next state $\tilde{Z} = (\tilde{Y},\tilde{\eb})$ is given by:
		\begin{align}
		\tilde{\yb}_v &= \yb_v + \hat{\yb}_v \label{eqNextState_y}\\
		\tilde{e}_v &= e_v + \hat{\yb}_v\T\1 + \sum_{u\in  \mathcal{N}_v} \big(\hat{e}_{ u\to v} - \hat{e}_{v \to u}\big), \label{eqNextState_e}
		\end{align}
		where the admissible action set of node $v$ is:
		\begin{align}
		\mathcal{A}_v(Z) = \mathcal{A}_v(\yb_v) := \{ \hat{\yb}\in \mathbb{R}^{|\mathcal{C}|}~|~\yb_v + \hat{\yb} \in \Omega_v  \}.
		\end{align}
		These dynamics satisfy:
		\begin{align}
		\sum_{v\in \mathcal{V}}\tilde{e}_v - \sum_{v\in \mathcal{V}}\tilde{\yb}_v\T\1= \sum_{v\in \mathcal{V}}{e}_v - \sum_{v\in \mathcal{V}}{\yb}_v\T\1.
		\end{align}
		
		\item[4.] \textit{Nodal cost function: } For a state $Z$ and an admissible action profile $\{\ab_v \in \mathcal{A}_v(\yb_v) \}_{v\in \mathcal{V}}$, the cost function of node $v$ is given by:
		\begin{equation}
		J_v(Z,\ab) = -\tilde{L}(\tilde{Y}) + \frac{\mu}{2} \sum_{u\in \mathcal{N}_v}[\tilde{e}_u]_+^2, \label{eqNodalCost}
		\end{equation}  
		where  $(\tilde{Y},\tilde{\eb})$ is the next state and $\mu>0$ is a penalty parameter. Here, the cost $J_v$ still involves the global (approximated) caching gain function $\tilde{L}$, but as we will show later, each node does not need to evaluate $J_v$. The second term in \eqref{eqNodalCost} represents a penalty of the differences in the estimation error terms between neighboring nodes. 
	\end{enumerate}

	This is a potential game with the potential function (to be minimized) given by
	\begin{equation}
	\Phi_{\mu}(Z,\ab) = -\tilde{L}(\tilde{Y}) + \frac{\mu}{2} \sum_{v\in \mathcal{V} } [\tilde{e}_v]_+^2. \label{eqPotentialFunc}
	\end{equation}
	This can be shown by noting that $\Phi_{\mu}(\tilde{Z},\0) = \Phi_{\mu}(Z,\ab)$ and that $\forall  \ab'_v \in \mathcal{A}_v(Z)$
	\begin{equation}
	\begin{split}
	&J_v(Z,\{\ab'_v, \ab_{-v}\}) - J_v(Z,\{\ab_v, \ab_{-v}\}) \\
	&\qquad = \Phi_{\mu}(Z,\{\ab'_v, \ab_{-v}\}) - \Phi_{\mu}(Z,\{\ab_v, \ab_{-v}\}),
	\end{split} 
	\label{eqMarginalPotential}
	\end{equation}
	where $\ab_{-v}$ denotes the actions of all the nodes other than $v$. 
	Condition \eqref{eqMarginalPotential} means that any improvement in the cost of node $v$ made by its local action is the same as the potential function improvement. 
	Moreover, it is easy to see that $\Phi_{\mu}$ is convex continuous with bounded level sets. 
	Thus, a stationary state Nash equilibrium always exists  and can be reached by the \emph{gradient play} strategy (see, e.g., \cite{marden2015game}). The following algorithm is an implementation of this strategy.

	\subsubsection{Algorithm Description} \label{secAlgorithm}
	We assume that time is partitioned into periods of equal length $T$, during which the nodes collect measurements from messages routed through them.
	Each node maintains and updates $\yb_v$ and the error term $e_v$ as follows.
	\begin{itemize}
		\item At period $t=0$, each node $v$ initializes $\yb_v(0) \in \Omega_v$ and $e_v(0) \gets (\1\T\yb_v(0) - c^0_v)$ such that 
		\begin{equation}
		\sum_{v\in \mathcal{V}}e_v(0) \le 0. \label{eqInitConstraint}
		\end{equation}
		\item At period $t>0$, node $v$ exchanges $e_v(t)$ with its neighbors and computes action $\ab_v(t)$:
		\begin{align}
		\hspace{-8mm} \hat{e}_{v\to u}(t) &= -\gamma_v(t)\frac{\partial J_v(Z(t),\ab)}{\partial \hat{e}_{v \to u}}\Big\rvert_{\ab=\0} \nnb\\
		&= \gamma_v(t)\mu \big([e_v(t)]_+ - [e_u(t)]_+ \big), \quad \forall u\in \mathcal{N}_v \label{eqAction_E}\\
		\hspace{-8mm} \hat{\yb}_v(t) &= \Big[ -\gamma_v(t) \nabla_{\hat{\yb}_{v}} J_v(Z(t),\ab) \big\rvert_{\ab=\0} \Big]_{\mathcal{A}_v(\yb_v(t))} \nnb\\
		&= \Big[  \gamma_v(t) \Big( \nabla_{\yb_{v}}\tilde{L}(Y)  -\mu\1[e_v(t)]_+  \Big) \Big]_{\mathcal{A}_v(\yb_v(t))} \label{eqAction_Y}
		\end{align}		
		where $\gamma_v(t)$ denotes the step size of node $v$ at iteration $t$. 
		Here, $\nabla_{\yb_{v}}\tilde{L}(Y)$ can be computed in a distributed fashion as shown in Section~\ref{DistributedComputationGrad} below. 
		Then node $v$ sends $\hat{e}_{v\to u}(t)$ to node $u\in \mathcal{N}_v$ and updates its state as follows: 
		\begin{align}
		\hspace{-8mm} \yb_v(t\!+\!1) &= \yb_v(t) + \hat{\yb}_v(t) \label{eqUpdate_Y}\\
		\hspace{-8mm}e_v(t\!+\!1) &= e_v(t) \!+\! \1\T \hat{\yb}_v(t) + \!\!\!\sum_{u\in \mathcal{N}_v}\! \big( \underbrace{\hat{e}_{ u \to v}(t)}_{\text{received}} - \underbrace{\hat{e}_{v \to u}(t)}_{\text{computed}} \big).
		\label{eqUpdate_E}
		\end{align}
	\end{itemize}
	It can be seen that the following forms are more convenient for implementation: 
	\begin{align}
	\yb_v(t+1) &=  \Big[\yb_v(t) + \gamma_v(t) \Big( \nabla_{\yb_{v}}\tilde{L}(Y)  -\mu [e_v(t)]_+\1  \Big) \Big]_{\Omega_v}  \label{eqUpdate_Y_proj}\\
	e_v(t+1) &= e_v(t) + \1\T(\yb_v(t+1) - \yb_v(t)) \nnb\\& +    \sum_{u\in \mathcal{N}_v}    \hat{e}_{ u \to v}(t) - \hat{e}_{v \to u}(t).  \label{eqUpdate_E_proj}
	\end{align}

	\begin{rem} \label{rem_decomposibility}
		In \cite{li2013designing,li2014decoupling,marden2015game}, the authors also provide a potential game-based algorithm for solving a (more general) constrained optimization problem, the design of which, if applied to \eqref{ProbCVX_Approx}, would yield an exponentially large state space. Specifically, to decompose $\Phi_{\mu}$, each node $v$ would need to keep track and update a local estimate $Y_v$ of the state $Y$ through exchanging information with direct neighbors. This would incur much more expensive communication and computational costs compared to our model and algorithm outlined above. Our advantage is gained by incorporating a distributed algorithm for each node to estimate partial gradients of $\tilde{L}$. Such an algorithm requires only a simple message exchange protocol, which we describe in Section~\ref{DistributedComputationGrad} below. 
	\end{rem}
	
	\subsection{Convergence}
	With initialization \eqref{eqInitConstraint} and updates \eqref{eqUpdate_Y}--\eqref{eqUpdate_E}, it can be shown that for any $t\ge 0$, 
	\begin{align*}
	\sum_{v\in \mathcal{V}}{e}_v(t) - \sum_{v\in \mathcal{V}} \1\T\yb_v(t) = \sum_{v\in \mathcal{V}}{e}_v(0) - \sum_{v\in \mathcal{V}} \1\T\yb_v(0) = M.
	\end{align*}
	Thus, 
	\begin{align}
	\sum_{v\in \mathcal{V}}{e}_v(t) = \sum_{v\in \mathcal{V}} \big( \1\T\yb_v(t) -c_v^0 \big), \quad \forall t \ge 0. \label{eq_sum_error}
	\end{align}
	The next result is obtained by following similar arguments as  \cite[Thm.~1]{li2014decoupling}. 	We provide details as our proof departs from  \cite{li2014decoupling} because our objective function is non-separable (while it is decoupled in \cite{li2014decoupling}). 
	
	\begin{theorem}\label{thm_Game_Equiv}
		For a fixed $\mu$, suppose a state action pair $\{Z,\ab\} = \{ (Y,\eb), (\hat{Y}, \hat{E}) \}$ is a stationary state Nash equilibrium. Then:
		\begin{itemize}
			\item[(i)] $Y$ is an optimizer of the following problem
			\begin{equation} \label{ProbCVX_Game_NE}
			\displaystyle{\max_{\{\yb_v \in \Omega_v\} }} \quad 
			\tilde{L}( Y)  - \frac{\mu}{2|\mathcal{V}|} \Big[ \sum_{v \in \mathcal{V}} (\yb_v\T \1 - c^0_v) \Big]_+^2.
			\end{equation}
			
			\item[(ii)] The estimation error $\eb$ satisfies
			$$
			[e_u]_+ = \frac{1}{|\mathcal{V}|}\Big[ \sum_{v \in \mathcal{V}} (\yb_v\T \1 - c^0_v) \Big]_+, \quad \forall u \in \mathcal{V}.
			$$
			\item[(iii)] The actions satisfy $\hat{\yb}_v = \0$ and $\sum_{u\in \mathcal{N}_v} \big(\hat{e}_{ u\to v} - \hat{e}_{v \to u}\big) = 0$ for all $v\in \mathcal{V}$. 
		\end{itemize}
	\end{theorem}
	\IEEEproof
	First, (iii) is obvious. Second, since  $\{Z,\ab\}$ is a stationary state Nash equilibrium, we have
	$$
	J_v(Z,\{\ab_v, \ab_{-v}\}) = \min_{\check{\ab}_v \in \mathcal{A}_v(Z)} J_v(Z,\{\check{\ab}_v, \ab_{-v}\}),\quad  \forall v\in \mathcal{V}.
	$$
	Since $J_v$ is convex and differentiable on $\check{\ab}_v = (\check{\yb}_v, \{\check{e}_{v\to u}\}_{u \in \mathcal{N}_v}) \in \mathcal{A}_v(Z) = \mathcal{A}_v(\yb_v)$, the condition above implies that for any $v\in \mathcal{V}$
	\begin{align*}
	\partial_{\check{e}_{v\to u} } J_v(Z,\{\check{\ab}_v, \ab_{-v}\})|_{\ab} &= 0, \\
	(\hat{\yb}'_v - \hat{\yb}_v)\T \Big[\nabla_{\check{\yb}_v} J_v(Z,\{\check{\ab}_v, \ab_{-v}\})|_{\ab}\Big] &\ge 0, \quad \forall \hat{\yb}'_v \in \mathcal{A}_v(\yb_v)
	\end{align*}
	which are respectively equivalent to
	\begin{align}
	&\!\!\![ \tilde{e}_v ]_+ - [\tilde{e}_u ]_+ = 0, \quad \forall u\in \mathcal{N}_v \label{eq_tilde_e_pos}\\
	&\!\!\!(\tilde{\yb}'_v - \tilde{\yb}_v)\T \big( \nabla_{\yb_{v}}\tilde{L}(\tilde{Y})  -\mu\1[\tilde{e}_v(t)]_+  \big) \ge 0, \quad \forall \tilde{\yb}'_v \in \Omega_v. \label{eqOptimalityCond_2}
	\end{align}
	Condition \eqref{eq_tilde_e_pos} and connectivity of $\mathcal{G}$ implies that $[ \tilde{e}_v ]_+ = [\tilde{e}_u ]_+$ for any $v,u\in \mathcal{V}$, which means either $\tilde{e}_v \le 0, \forall v\in \mathcal{V}$ or $\tilde{e}_v = \tilde{e}_u \ge 0, \forall v,u\in \mathcal{V}$. In any case, the following holds for any $u\in \mathcal{V}$ 
	$$
	[ \tilde{e}_u ]_+ = \frac{1}{|\mathcal{V}|} \Big[\sum_{v \in \mathcal{V}}\tilde{e}_v \Big]_+ \stackrel{\eqref{eq_sum_error}}{=}\frac{1}{|\mathcal{V}|} \Big[ \sum_{v \in \mathcal{V}} \big(\tilde{\yb}_v\T\1 -c_v^0 \big) \Big]_+.
	$$
	This together with (iii) proves (ii). It remain to show (i). From~\eqref{eqOptimalityCond_2} and using (ii) and (iii), we have for all $v\in \mathcal{V}$ and $\forall \tilde{\yb}'_v \in \Omega_v$
	\begin{align*}
	(\tilde{\yb}'_v - {\yb}_v)\T \Big( \nabla_{\yb_{v}}\tilde{L}(Y)  -\frac{\mu}{|\mathcal{V}|} \1 \Big[\sum_{v \in \mathcal{V}} ({\yb}_v\T\1 -c_v^0) \Big]_+  \Big) &\ge 0.
	\end{align*}
	This clearly shows that ${\yb}_v$ is optimal to \eqref{ProbCVX_Game_NE}. 
	\endIEEEproof
	
	The following result is obvious from Theorem~\ref{thm_Game_Equiv}-(i). 
	\begin{corollary}\label{corInftyPenalty}
		As $\mu \to \infty$, the equilibria of the game constitute solutions to \eqref{ProbCVX_Game}--\eqref{eqConstraintTotalCacheSize}. 
	\end{corollary}
	
	Before proving the convergence of the algorithm, we summarize approximation steps introduced so far in dealing with the original problem (MaxCG). First, we relax the binary constraints \eqref{BinaryConstraint} and approximate the objective function $F$ by $L$ in \eqref{eqL}, thereby obtaining \eqref{ProbCVX}, a convex problem on the relaxed feasible set. Second, since $L$ is nondifferentiable, we then replace it with $\tilde{L}$ in \eqref{eqL3}.  Third, by resorting to the potential game theory, we effectively remove the global constraint \eqref{GlobalCap} by adding a penalizing term to $\tilde{L}$, resulting \eqref{ProbCVX_Game_NE}. In summary, we have the following approximations in terms of caching gains. 
	\begin{equation}\label{eq3StepApprox}
	(\textrm{MaxCG}) ~~\approx~~ \eqref{ProbCVX} ~~\approx~~ \eqref{ProbCVX_Game}-\eqref{eqConstraintTotalCacheSize}
	~~\approx~~ \eqref{ProbCVX_Game_NE}.
	\end{equation}
	
	\begin{theorem}
		For any small $\epsilon>0$, there exist $\mu$ sufficiently large and $\alpha$ sufficiently small such that \eqref{ProbCVX_Game_NE} approximates (MaxCG) within $(1-1/e-\epsilon)$-ratio in terms of the caching gain. 
	\end{theorem}
	\IEEEproof
	First, the approximation ratio of the first step in \eqref{eq3StepApprox} is $(1-1/e)$; see \eqref{eqBound_F}. Second, the approximation errors in  the last two steps can be made arbitrarily small by choosing sufficiently small $\alpha$ and large $\mu$; see \eqref{eqL3Bound} and Corollary~\ref{corInftyPenalty}. Thus, we conclude that the 3-step approximation in \eqref{eq3StepApprox} can achieve ratio $(1-\epsilon-1/e)$. 
	\endIEEEproof
	
	The following result establishes the convergence of the above algorithm for a uniform constant step size. 
	
	\begin{theorem} 
		Consider the algorithm described in Section~\ref{secAlgorithm} with $\gamma_v(t) \equiv \gamma$ such that 
		\begin{equation}
		\gamma < \bar{\gamma}^0 := \frac{2}{\alpha^{-1}C_0 + 2\mu}.
		\label{eqLambda_0}
		\end{equation}
		Then, $(Z(t),\ab(t))$ converges to a stationary state Nash equilibrium. Moreover, any limit point $Y^*$ of $\{Y(t)\}$ is an optimizer of \eqref{ProbCVX_Game_NE}. 
	\end{theorem}
	\IEEEproof
	Note that the gradient $\nabla_{\ab} \Phi_{\mu}(Z,\ab)$ is Lipschitz continuous with a parameter
	\begin{equation}
	K_{\nabla \Phi} = K_{\nabla \tilde{L}} + 2\mu =  \alpha^{-1}C_0 + 2\mu.
	\end{equation}
	Thus, it follows from 
	\cite[Thm.~4]{li2013designing}  
	that the algorithm converges to a stationary state Nash equilibrium 
	for any $\gamma \in (0, \bar{\gamma}^0)$. The rest of the proof follows from Theorem~\ref{thm_Game_Equiv}. 
	\endIEEEproof
	
	Note that $\bar{\gamma}^0$ is a theoretical bound for the gradient method, while step sizes larger than $\bar{\gamma}^0$ often still work in practice; of course, the larger the step sizes are, the closer to instability the algorithm is.
	In this paper, we focus on the case of uniform step size and synchronous communications, but it can be shown further that the algorithm is also robust to bounded communication delays, asynchronism of the nodes' clocks (or update times), and heterogeneous and time-varying nodal step sizes; see, e.g., \cite[Chap.~3 and 7]{Bertsekas1989parallel} and \cite{marden2012state}.

	\begin{rem} \label{rem_constraint_violation} 
		Given a fixed $\mu$, the global capacity constraint in \eqref{eqConstraintTotalCacheSize} is likely to be violated due to the penalizing term in \eqref{ProbCVX_Game_NE}.
		To reduce such violation, we can initialize $\sum_{v\in \mathcal{V}}c_v^0 = (M-\epsilon)$ for some small $\epsilon \in (0,1)$ and select $\mu = O(\tilde{L}^*)$. Since $\tilde{L}^*$ is unknown, we can choose $\mu = O(C_0)$ (noting that $L({Y}) \le C_0, \forall Y$), where $C_0$ (or an upper bound  $\bar{C}_0$) can be estimated in a centralized fashion from history data or in a distributed manner as described in Section~\ref{secParameterEstimate} below.
	\end{rem}

	\subsection{Implementation Considerations}
	This subsection details on how each node in the network can obtain information needed to implement the algorithm. This includes: online estimations of partial derivatives $\partial_{y_{vi}} \tilde{L}$, step size bound in \eqref{eqLambda_0} (for ensuring convergence), and an eviction policy for updating cache contents.
	
	\subsubsection{Distributed gradient estimation}  \label{DistributedComputationGrad}
	We adopt the mechanism used in \cite{ioannidis2016adaptive, gill2016bidcache}, namely, additional control messages are attached to the request and response traffic to gather needed information. This enables each node $v$ to estimate partials $\partial_{\yb_{v}} \tilde{L}$ in a distributed fashion by using information in the messages passing by during each time interval $T$. In particular:
	\begin{itemize}
		\item Every time a node generates a new request $(i,p) \in \mathcal{R}$, it also creates an additional control message $m_s$ to send over $p$ along with the request. At node $p_1$, $m_s(p_1) = y_{p_1 i}$. As this message is propagated to node $p_l$, $m_s$ is updated as follows:
		\begin{equation}\label{MsgMs}
		m_s(p_l) = m_s(p_{l-1}) + y_{p_l i}
		\end{equation}
		until a node $u\in p$ such that $m_s(u) > 1+\frac{\alpha}{2}$ is found or the end of the path is reached (in which case $u=p_{|p|}$). Each visited node $p_l$  keeps a local copy of $m_s(p_l)$. 
		
		\item Node $u$ (found above) creates a control message $m_r$ to send back in the reverse direction. At $u$, $m_r(u)  = 0$. At $p_l$, 
		\begin{equation}\label{MsgMr}
		m_r(p_l) = m_r(p_{l+1}) + w_{p_{l+1}p_{l}}\mathrm{sat}'_{\alpha}(m_s(p_l)),
		\end{equation}
		where $\mathrm{sat}'_{\alpha}(x) = \frac{d}{dx}\mathrm{sat}_{\alpha}(x)$, i.e., 
		$$
		\mathrm{sat}'_{\alpha}(x) = \begin{cases}
		0 & \text{if } x\ge 1+\frac{\alpha}{2}\\
		\frac{1}{\alpha}(1+\frac{\alpha}{2}-x) & \text{if } 1-\frac{\alpha}{2} \le x <1+\frac{\alpha}{2}\\
		1 & \text{if } 0\le x< 1-\frac{\alpha}{2}.
		\end{cases}
		$$
		
		\item For each item $i$ and each node $v$, let $$t_{vi} := m_r(v)$$ as computed above.  It can be seen that $t_{vi}$ is proportional to the partial derivative of $\tilde{L}$ for request $(i,p) \in \mathcal{R}$, i.e., 
		$$
		t_{vi} = \frac{\partial}{\partial y_{vi}} \sum_{k=k_p(v)}^{|p|-1} w_{p_{k+1} p_k} \mathrm{sat}_{\alpha} \big(\sum_{l=1}^k y_{p_{l}i} \big),
		$$
		where $k_p(v)$ denotes the position of $v$ in $p$.
		
		\item It remains to show how each node $v$ estimates the partial derivative of $\tilde{L}$ with respect to $y_{vi}, \forall i\in \mathcal{C}$.
		This is trivial if the rate $\lambda_{(i,p)}$ is known to all the nodes in path $p$; otherwise, each node needs to estimate it. 
		To this end, let $\mathcal{T}_{vi}$ denote the set of $t_{vi}$ collected by node $v$ regarding item $i$ during each time slot. Then it can be shown \cite[Lem.~1]{ioannidis2016adaptive} that 
		$$
		z_{vi} := \frac{\sum_{t\in \mathcal{T}_{vi}}t}{T}
		$$
		is an unbiased estimate of the partial derivative $\partial_{y_{vi}} \tilde{L}$.  
	\end{itemize}

	\subsubsection{Distributed estimation of $\gamma_0$} \label{secParameterEstimate}
	To implement the algorithm, all nodes need to agree not only on a common $\mu$ and $\alpha$, but also on the step size bound $\bar{\gamma}_0$ in \eqref{eqLambda_0}; the latter depends on $C_0$, where we recall that $C_0 = \sum_{(i,p)\in \mathcal{R}} \lambda_{(i,p)} \sum_{k=1}^{|p|-1} w_{p_{k+1} p_k}$. Thus, we now focus on how to estimate $C_0$ or an upper bound in a distributed fashion.
	
	First, we assume that each node $v\in \mathcal{V}$ knows the weight of the path $$w_{p} := \sum_{k=1}^{|p|-1} w_{p_{k+1} p_k}$$ and an estimate (or an upper bound) of the associated rate, denoted by $\bar{\lambda}_{(i,p)}$, for any $(i,p) \in \mathcal{R}$ such that $p_1 = v$. In fact, node $v$ can compute $w_p$ simply by probing path $p$ and the end node replying with a control message sent in the reverse direction to accumulate the weight of the path (this can be done a priori or periodically).
	
	Then, every node $v\in \mathcal{V}$ can find
	$$
	C_{v0} := \sum_{(i,p)\in \mathcal{R}, p_1=v} \bar{\lambda}_{(i,p)}w_{p}.
	$$
	Clearly, $\sum_{v\in \mathcal{V}}C_{v0} \ge C_0$. Moreover, by running an additional average consensus algorithm (e.g., \cite{xiao2004fast,olfati2007consensus}) with initial conditions $\{C_{v0} \}_{v\in \mathcal{V}}$, all the nodes in $\mathcal{V}$ can compute $\frac{\sum_{v\in \mathcal{V}}C_{v0} }{|\mathcal{V}|}$. For completeness, a detailed algorithm is given in Appendix~\ref{secConsensusAlg}. It should be noted that such algorithm converges exponentially fast and independently from our main algorithm described above. 
	Therefore, assuming that an upper bound $\bar{N}$ on $|\mathcal{V}|$ is known to all the nodes, they can find
	\begin{equation}
	\bar{C}_0 := \bar{N} \frac{\sum_{v\in \mathcal{V}}C_{v0} }{|\mathcal{V}|} \label{eqC_0_bar}
	\end{equation} which clearly satisfies $\bar{C}_0 \ge C_0$.
	
	Finally, we will choose $\mu = \mu_0\bar{C}_0$ for some $\mu_0$ chosen a priori together with $\alpha$, thereby having
	\begin{equation}
	\tilde{\gamma}^0 := \frac{2}{(\alpha^{-1} + 2\mu_0)\bar{C}_0} \le \bar{\gamma}^0.
	\label{eqLambda_0_tilde}
	\end{equation}
	
	\subsubsection{Eviction policy}\label{subsecEviction}
	At the end of each iteration $t$, before deciding what to put in the cache, each node $v$ needs to determine the maximum number of items that it can store. This number is based on the expected local caching capacity $\sum_{i\in \mathcal{C}} y_{vi}(t)$, which can be fractional. 
	Therefore, a local rounding scheme is needed; e.g., randomized rounding as in \cite{ioannidis2016adaptive}.  
	A simpler heuristic would be the following: node $v$ determines a positive integer $c_v(t) \in [0,\bar{c}_v]$ such that:
	\begin{itemize}
		\item If the global constraint \eqref{GlobalCap} is a hard constraint, then
		$$ c_v(t) = \min\{ \lfloor\1\T\yb_v(t) \rfloor, \bar{c}_v\}.$$
		\item If \eqref{GlobalCap} is a soft constraint, then $c_v(t)$ is the nearest to the sum $\sum_{i\in \mathcal{C}} y_{vi}(t)$, i.e., 
		$$ c_v(t) = \min\big\{\big[\1\T\yb_v(t)\big]_\mathbb{N}, \bar{c}_v \big\}.$$
	\end{itemize}
	Node $v$ then places at most $c_v(t)$ content items, corresponding to the largest elements of $\yb_v(t)$, into its cache.
	
	\subsubsection{Efficient transmission of control messages}
	Our algorithm requires each node to perform only few basic operations at each iteration to update its states \eqref{eqUpdate_Y_proj}--\eqref{eqUpdate_E_proj} and traversing control messages \eqref{MsgMs}--\eqref{MsgMr} for estimating local gradients; the projection $[\cdot]_{\Omega_v}$ in \eqref{eqUpdate_Y_proj} can be as simple as scaling. 
	We note that the control messages $\hat{e}_{v\to u}$, $m_s$, and $m_r$, described in \eqref{eqAction_E}, \eqref{MsgMs}, and \eqref{MsgMr} respectively, can be encoded in very few bytes and in many cases can be piggybacked onto the existing traffic of Interest and Data packets that normally flow through each node.
	In particular, the message $m_s$ generated for request $(i,p) \in \mathcal{R}$ will be attached to the Interest packet for item $i$ for as many hops as possible, while the corresponding reply $m_r$ will be attached, if possible, to the Data packet containing item $i$. Message $\hat{e}_{v\to u}, u\in \mathcal{N}_v$ can be attached to any packet being transmitted from $v$ to $u$, regardless of its type. These messages do not need to be transmitted immediately if the link between $v$ and $u$ is idle: in this case they can be placed in a queue where they wait until the next available transmit opportunity or until a timeout expires, whichever occurs first. As a consequence, the overall overhead and storage of these messages are negligible.
	
	In the case of an NDN network, since these messages are propagated in a hop-by-hop fashion, we recommend encoding them as NDNLPv2 header fields. If necessary, multiple messages can be attached to the same NDNLPv2 packet, providing further bandwidth savings.
	
	We also note that the loss of one or more control messages can reduce the convergence rate of the algorithm, but will not affect its correctness.

	\section{Numerical Examples}
	\label{sec:sim}
	
	In this section, we demonstrate the performance of our algorithm applied to several network topologies.
	
	\paragraph{Topologies: }
	We consider the networks shown in Table~\ref{tbl_topologies}. \verb|grid_2d| is a two-dimensional square grid and \verb|expander| is a Margulies-Gabber-Galil expander~\cite{gabber1981explicit}. The next four graphs are random graphs sampled from a distribution. 
	\verb|erdos_renyi| is an Erdos-Renyi graph with parameter $p=0.1$; 
	\verb|small_world| is a small-world graph~\cite{kleinberg2000small} that consists of a grid with additional long-range links; graph 
	\verb|watts_strogatz| is generated according  to the Watts-Strogatz model in~\cite{watts1998collective}; and 
	\verb|barabasi_albert| follows the preferential attachment model in~\cite{barabasi1999emergence}. 
	The last three graphs are the GEANT, Abilene, and Deutsche Telekom backbone networks~\cite{rossi2011caching}.
	
	\begin{table}
		\caption{Graph topologies and parameters}
		\label{tbl_topologies}
		\centering
		\begin{tabular}{lllllll}
			\rule[-1ex]{0pt}{2.5ex}  Graph & $|\mathcal{V}|$ & $|\mathcal{E}|$ & $|\mathcal{C}|$ & $|\mathcal{Q}|$ & $|\mathcal{R}|$ & $M$ \\ 
			\hline 
			\rule[-1ex]{0pt}{2.5ex}  \verb|grid_2d (G2)|			& 100 & 180 & 100 & 20 & 1K & 300  \\ 
			\rule[-1ex]{0pt}{2.5ex}  \verb|expander (Ex)|			& 100 & 340 & 100 & 50 & 2K & 400   \\ 
			\hline 
			\rule[-1ex]{0pt}{2.5ex}  \verb|barabasi_albert (BA)|	& 100 & 384 & 100 & 50 & 2K & 400   \\ 
			\rule[-1ex]{0pt}{2.5ex}  \verb|small_world (SW)|		& 100 & 240 & 100 & 50 & 2K & 400   \\ 
			\rule[-1ex]{0pt}{2.5ex}  \verb|watts_strogatz (WS)|		& 100 & 200 & 100 & 50 & 2K & 400   \\ 
			\rule[-1ex]{0pt}{2.5ex}  \verb|erdos_renyi (ER)|		& 100 & 521 & 100 & 50 & 2K & 400   \\ 
			\hline 
			\rule[-1ex]{0pt}{2.5ex}  \verb|geant (Ge)|				& 22  & 33  & 100 & 20 & 1K & 144    \\ 
			\rule[-1ex]{0pt}{2.5ex}  \verb|abilene (Ab)|			& 9   & 13  & 10  & 9  & 100& 28    \\ 
			\rule[-1ex]{0pt}{2.5ex}  \verb|dtelekom (Dt)|			& 68  & 273 & 100 & 20 & 1K & 304   \\ 
		\end{tabular}%
	\end{table}

	\paragraph{Experiment setup: }
	For each graph, we generate a catalog $\mathcal{C}$ and assign each item $i\in \mathcal{C}$ to a node selected uniformly at random (u.a.r.) from $\mathcal{V}$. We select the weight of each edge u.a.r. from $[0.01,1]$ and a set of consumers $\mathcal{Q}\subset \mathcal{V}$ u.a.r. Each consumer $v \in \mathcal{Q}$ requests an item $i$ selected from $\mathcal{C}$ according to a Zipf distribution with parameter $1.2$. The request is routed over the shortest path $p$ between $v$ and the designated server for item $i$. The set of requests is denoted by $\mathcal{R}$. 
	We choose $\bar{c}_v = |\mathcal{C}|$ and measurement/update period $T=1$. 
	Moreover, we let:
	\begin{itemize}
		\item $\alpha = 0.2$. Thus, it follows from \eqref{eqL3Bound} that $|\tilde{L}({Y}) - L({Y})| \le 2.5\% C_0, \forall Y\in \mathcal{D}_2$. The actual error is often much less.
		\item $\mu = \displaystyle\frac{1}{4}\bar{C}_0$. Thus, by \eqref{eqLambda_0_tilde}, $\tilde{\gamma}^0 = \displaystyle\frac{4}{11\bar{C}_0}.$ 
		\item step size $\gamma_v \equiv \tilde{\gamma}_0, \forall v\in \mathcal{V}$. 
		\item $c_v^0 = \displaystyle\frac{M-\epsilon}{|\mathcal{V}|}$, with $\epsilon = 0.1$, see Remark~\ref{rem_constraint_violation}.  
	\end{itemize}

\begin{figure}
		\centering
		\includegraphics[width=0.95\columnwidth,trim=25 0 5 10,clip]{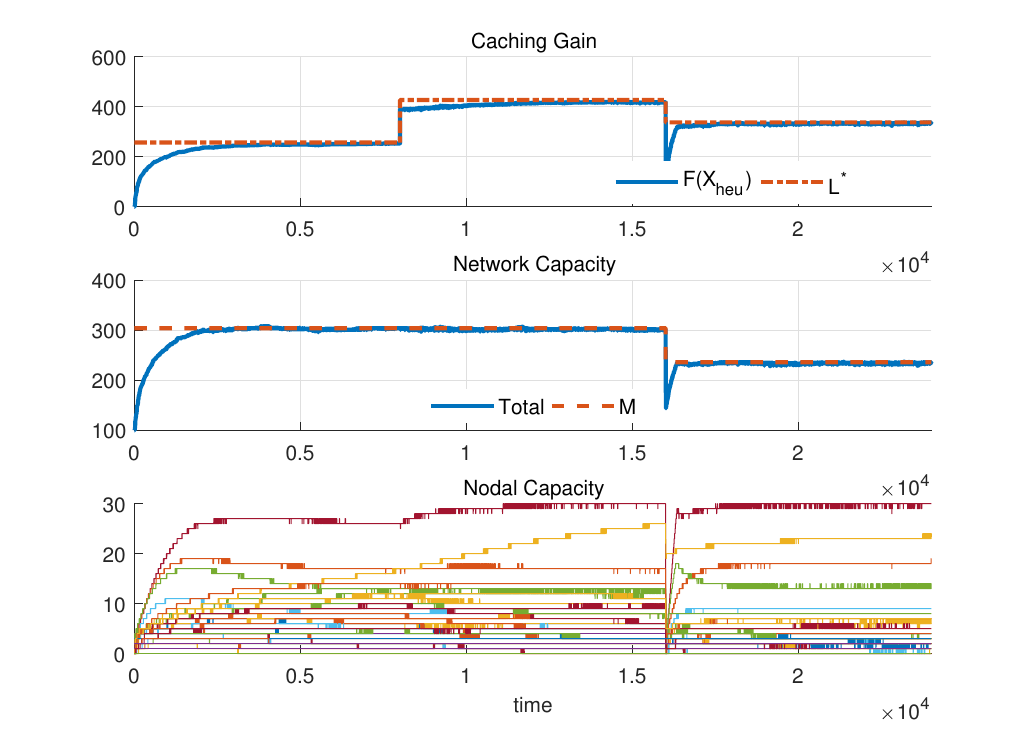}
		\caption{(Color online) Simulation results for $\sf{dtelekom}$ network using cache allocation $X_{\textrm{heu}}$ obtained from our heuristic placement in Sect.~\ref{subsecEviction}. $L^*$ from~\eqref{ProbCVX} is obtained by a centralized algorithm and is an upper bound on the optimal caching gain. Bottom plot shows cache sizes $c_v(t), \forall v \in \mathcal{V}$.}
		\label{fig_dtelekom}
	\end{figure}

	\paragraph{Results: }
	First, we simulate our algorithm on the \verb|dtelekom| graph. During time interval $[0,8000]$, the request rates $\lambda_{(i,p)}$ are selected u.a.r. from $[0.1,1]$ and after that $\lambda_{(i,p)}=1,  \forall (i,p) \in \mathcal{R}$. We will use $\bar{\lambda}_{(i,p)}=1$ as an upper bound of $\lambda_{(i,p)}, \forall (i,p) \in \mathcal{R}$ for computing $\bar{C}_0$ as in \eqref{eqC_0_bar} (assuming $\bar{N} = |\mathcal{V}|$). Moreover, we reduce the budget from $M$ to $(M-|\mathcal{V}|)$ at $t=16000$. The simulation results are shown in Figure~\ref{fig_dtelekom}, which clearly demonstrates adaptability and optimality of our algorithm.
	
	As we can observe, from initial allocation $Y(0) = X(0)$, the network quickly reaches total cache size $M$. After that, the caching gain is improved  and nearly reaches upper bound ${L^*}$, thereby implying near optimality. 
	
	\begin{rem} (\emph{On adaptability}) 
		From simulations, the convergence rate of our algorithm seems to be sublinear (expected since $\tilde{L}$ is not strongly concave). Thus, our algorithm is suitable for networks with not too fast changes. 
	\end{rem}

\begin{figure}[!tb]
		\centering
		\includegraphics[width=0.95\columnwidth,trim=34 0 25 0,clip]{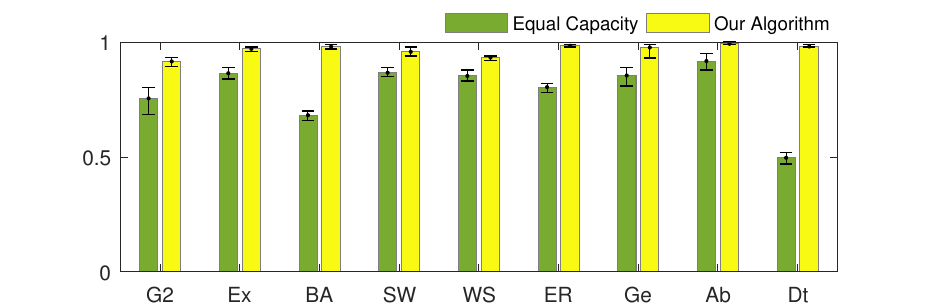}
		\caption{Comparison of normalized caching gains for graphs in Table~\ref{tbl_topologies}. We show here the average and error bars of 10 runs for each scenario.}
		\label{fig_compare_ratio}
	\end{figure}

	We also compare the performance, in terms of caching gains (normalized to $L^*$), of our algorithm with the centralized solution approach using the equal node-capacity allocation across all topologies in Table~\ref{tbl_topologies}. 
	Specifically, the latter fixes  
	$${c}_v(t)  \equiv  \bar{c}_v  =  \frac{M-|\mathcal{C}|}{|\mathcal{V}|}  +  |\{i: v \in  \mathcal{S}_i\}|, \quad\forall v \in  \mathcal{V},$$ i.e., \eqref{GlobalCap} is redundant as $\sum_{v\in\mathcal{V}}\bar{c}_v  =  M$. 
	Note that the optimal (relaxed) caching gain in equal node capacity, denoted by $L_{EC}^*$ and obtained by solving \eqref{ProbCVX}  without global constraint~\eqref{GlobalCap}, is not only an upper bound on caching gains of all  suboptimal caching policies in the same setting, but also a lower bound of $L^*$ in~\eqref{ProbCVX} with global constraint~\eqref{GlobalCap} and $\bar{c}_v = |\mathcal{C}|$. Significant gaps (ranging from $15\%$ to $50\%$) between $L^*_{EC}$ and other common caching strategies  have been shown in \cite{ioannidis2016adaptive} for a similar set of topologies. Here, we focus on showing improvement of $F(X_{heu})$ over $L_{EC}^*$. To this end, we run our algorithm  for $10^4$ time units with $\lambda_{(i,p)}  = 1, \forall (i,p)  \in  \mathcal{R}$ and estimate the steady state caching gain by averaging  the objective values $F(X_{heu})$ over the last $10^3$ time units. Fig.~\ref{fig_compare_ratio} shows the average results of 10 runs, which clearly demonstrate that our algorithm yields (near) optimal caching gains and outperform the best centralized solutions with equal capacity across all the topologies considered.
	
	\section{Conclusions} \label{sec:conclusion}
    
    We have designed a distributed and adaptive ICN caching scheme with optimality guarantees. Previous work in this area assumes that per-node cache sizes are predetermined constants and focuses on content placement. Our novel contribution addresses the problem of dynamic cache size design in emerging cloud-based networks that maximizes performance under varying network conditions. The resulting decentralized algorithm converges to cache allocations that are within a factor $1-1/e-\epsilon$ from the optimal. In addition to optimal content placement, the algorithm reallocates a given network-wide cache budget among the nodes as needed to maintain optimal cache allocation as network user demand changes. While we assumed that requests for any content item have a predetermined (typically shortest) path to a producer, enhancing our solution to take into account dynamic non-shortest path routing (e.g., to avoid congested paths) is an important open question.

    
    \bibliographystyle{IEEEtran}

\begin{thebibliography}{10}
    	\providecommand{\url}[1]{#1}
    	\csname url@rmstyle\endcsname
    	\providecommand{\newblock}{\relax}
    	\providecommand{\bibinfo}[2]{#2}
    	\providecommand\BIBentrySTDinterwordspacing{\spaceskip=0pt\relax}
    	\providecommand\BIBentryALTinterwordstretchfactor{4}
    	\providecommand\BIBentryALTinterwordspacing{\spaceskip=\fontdimen2\font plus
    		\BIBentryALTinterwordstretchfactor\fontdimen3\font minus
    		\fontdimen4\font\relax}
    	\providecommand\BIBforeignlanguage[2]{{%
    			\expandafter\ifx\csname l@#1\endcsname\relax
    			\typeout{** WARNING: IEEEtran.bst: No hyphenation pattern has been}%
    			\typeout{** loaded for the language `#1'. Using the pattern for}%
    			\typeout{** the default language instead.}%
    			\else
    			\language=\csname l@#1\endcsname
    			\fi
    			#2}}
    	
    	\bibitem{carofiglio2015lac}
    	G.~Carofiglio, L.~Mekinda, and L.~Muscariello, ``{LAC}: Introducing
    	latency-aware caching in information-centric networks,'' in \emph{Proc. 40th
    		Conf. Loc. Computer Netw.}\hskip 1em plus 0.5em minus 0.4em\relax IEEE, 2015,
    	pp. 422--425.
    	
    	\bibitem{thomas2015object}
    	Y.~Thomas, G.~Xylomenos, C.~Tsilopoulos, and G.~C. Polyzos, ``Object-oriented
    	packet caching for {ICN},'' in \emph{Proc. 2nd ACM Conf. Info.-Centric
    		Networking}.\hskip 1em plus 0.5em minus 0.4em\relax ACM, 2015, pp. 89--98.
    	
    	\bibitem{nguyen2015congestion}
    	D.~Nguyen, K.~Sugiyama, and A.~Tagami, ``Congestion price for cache management
    	in information-centric networking,'' in \emph{Proc. IEEE Conf. Computer
    		Commun. Wkshps}.\hskip 1em plus 0.5em minus 0.4em\relax IEEE, 2015, pp.
    	287--292.
    	
    	\bibitem{chai2012cache}
    	W.~K. Chai, D.~He, I.~Psaras, and G.~Pavlou, ``Cache less for more  in
    	information-centric networks,'' in \emph{Proc. Int. Conf. Research
    		Networking}.\hskip 1em plus 0.5em minus 0.4em\relax Springer, 2012, pp.
    	27--40.
    	
    	\bibitem{dehghan2016utility}
    	M.~Dehghan, L.~Massoulie, D.~Towsley, D.~Menasche, and Y.~C. Tay, ``A utility
    	optimization approach to network cache design,'' in \emph{Proc. 35th Annu.
    		IEEE Int. Conf. Computer Commun.}, 2016, pp. 1--9.
    	
    	\bibitem{ming2014age}
    	Z.~Ming, M.~Xu, and D.~Wang, ``Age-based cooperative caching in
    	information-centric networking,'' in \emph{23rd Int. Conf. Computer Commun.
    		Netw.}, 2014, pp. 1--8.
    	
    	\bibitem{badov2014congestion}
    	M.~Badov, A.~Seetharam, J.~Kurose, V.~Firoiu, and S.~Nanda, ``Congestion-aware
    	caching and search in information-centric networks,'' in \emph{Proc. 1st ACM
    		Conf. Info.-Centric Network.}, 2014, pp. 37--46.
    	
    	\bibitem{ioannidis2016adaptive}
    	S.~Ioannidis and E.~Yeh, ``Adaptive caching networks with optimality
    	guarantees,'' in \emph{ACM SIGMETRICS Performance Evaluation Rev.}, vol.~44,
    	no.~1, 2016, pp. 113--124.
    	
    	\bibitem{che2002hierarchical}
    	H.~Che, Y.~Tung, and Z.~Wang, ``Hierarchical web caching systems: Modeling,
    	design and experimental results,'' \emph{Selected Areas in Communications},
    	vol.~20, no.~7, pp. 1305--1314, 2002.
    	
    	\bibitem{jacobson2009networking}
    	V.~Jacobson, D.~K. Smetters, J.~D. Thornton, M.~F. Plass, N.~H. Briggs, and
    	R.~L. Braynard, ``Networking named content,'' in \emph{Proc. 5th Int. Conf.
    		Emerging Networking Experim. Tech.}, 2009, pp. 1--12.
    	
    	\bibitem{lv2002search}
    	Q.~Lv, P.~Cao, E.~Cohen, K.~Li, and S.~Shenker, ``Search and replication in
    	unstructured peer-to-peer networks,'' in \emph{Proc. 16th Int. Conf.
    		Supercomputing}, 2002, pp. 84--95.
    	
    	\bibitem{rossi2011caching}
    	D.~Rossi and G.~Rossini, ``Caching performance of content centric networks
    	under multi-path routing (and more),'' Telecom ParisTech, Tech. Rep., 2011.
    	
    	\bibitem{krause2012}
    	A.~Krause and D.~Golovin, ``Submodular function maximization,''
    	\emph{Tractability: Practical Approaches to Hard Problems}, vol.~3, no.~19,
    	p.~8, 2012.
    	
    	\bibitem{greedy2}
    	G.~L. Nemhauser, L.~A. Wolsey, and M.~L. Fisher, ``An analysis of
    	approximations for maximizing submodular set functions---i,''
    	\emph{Mathematical Programming}, vol.~14, no.~1, pp. 265--294, Dec 1978.
    	
    	\bibitem{vondrak2008optimal}
    	J.~Vondr{\'a}k, ``Optimal approximation for the submodular welfare problem in
    	the value oracle model,'' in \emph{STOC}, 2008.
    	
    	\bibitem{calinescu2007maximizing}
    	G.~Calinescu, C.~Chekuri, M.~P{\'a}l, and J.~Vondr{\'a}k, ``Maximizing a
    	submodular set function subject to a matroid constraint,'' in \emph{Integer
    		programming and combinatorial optimization}.\hskip 1em plus 0.5em minus
    	0.4em\relax Springer, 2007, pp. 182--196.
    	
    	\bibitem{calinescu2011}
    	------, ``Maximizing a monotone submodular function subject to a matroid
    	constraint,'' \emph{SIAM Journal on Computing}, vol.~40, no.~6, pp.
    	1740--1766, 2011.
    	
    	\bibitem{nemhauser1978best}
    	G.~L. Nemhauser and L.~A. Wolsey, ``Best algorithms for approximating the
    	maximum of a submodular set function,'' \emph{Mathematics of operations
    		research}, vol.~3, no.~3, pp. 177--188, 1978.
    	
    	\bibitem{shanmugam2013femtocaching}
    	K.~Shanmugam, N.~Golrezaei, A.~G. Dimakis, A.~F. Molisch, and G.~Caire,
    	``Femtocaching: Wireless content delivery through distributed caching
    	helpers,'' \emph{Transactions on Information Theory}, vol.~59, no.~12, pp.
    	8402--8413, 2013.
    	
    	\bibitem{ioannidis2018adaptive}
    	S.~Ioannidis and E.~Yeh, ``Adaptive caching networks with optimality
    	guarantees,'' in \emph{Transactions on Networking}, 2018.
    	
    	\bibitem{ageev2004pipage}
    	A.~A. Ageev and M.~I. Sviridenko, ``Pipage rounding: A new method of
    	constructing algorithms with proven performance guarantee,'' \emph{Journal of
    		Combinatorial Optimization}, vol.~8, no.~3, pp. 307--328, 2004.
    	
    	\bibitem{ioannidis2017icn}
    	S.~Ioannidis and E.~Yeh, ``Jointly optimal routing and caching for arbitrary
    	network topologies,'' in \emph{ACM ICN}, 2017.
    	
    	\bibitem{ioannidis2018jointly}
    	------, ``Jointly optimal routing and caching for arbitrary network
    	topologies,'' \emph{IEEE Journal on Selected Areas in Communications, Special
    		Issue on Caching for Communications and Networks}, 2018.
    	
    	\bibitem{palomar2006tutorial}
    	D.~P. Palomar and M.~Chiang, ``A tutorial on decomposition methods for network
    	utility maximization,'' \emph{IEEE J. Sel. Areas. Commun.}, vol.~24, no.~8,
    	pp. 1439--1451, 2006.
    	
    	\bibitem{marden2012state}
    	J.~R. Marden, ``State based potential games,'' \emph{Automatica}, vol.~48,
    	no.~12, pp. 3075--3088, 2012.
    	
    	\bibitem{li2013designing}
    	N.~Li and J.~R. Marden, ``Designing games for distributed optimization,''
    	\emph{{IEEE} J. Sel. Topics Signal Process.}, vol.~7, no.~2, pp. 230--242,
    	2013.
    	
    	\bibitem{li2014decoupling}
    	------, ``Decoupling coupled constraints through utility design,'' \emph{{IEEE}
    		Trans. Autom. Control}, vol.~59, no.~8, pp. 2289--2294, 2014.
    	
    	\bibitem{marden2015game}
    	J.~R. Marden and J.~S. Shamma, ``Game theory and distributed control,'' in
    	\emph{Handbook of game theory with economic applications}.\hskip 1em plus
    	0.5em minus 0.4em\relax Elsevier, 2015, vol.~4, pp. 861--899.
    	
    	\bibitem{li2011virtual}
    	W.~Li, J.~Tordsson, and E.~Elmroth, ``Virtual machine placement for predictable
    	and time-constrained peak loads,'' in \emph{International Workshop on Grid
    		Economics and Business Models}.\hskip 1em plus 0.5em minus 0.4em\relax
    	Springer, 2011, pp. 120--134.
    	
    	\bibitem{guenter2011managing}
    	B.~Guenter, N.~Jain, and C.~Williams, ``Managing cost, performance, and
    	reliability tradeoffs for energy-aware server provisioning,'' in
    	\emph{INFOCOM, 2011 Proceedings IEEE}.\hskip 1em plus 0.5em minus 0.4em\relax
    	IEEE, 2011, pp. 1332--1340.
    	
    	\bibitem{van2010cost}
    	R.~Van~den Bossche, K.~Vanmechelen, and J.~Broeckhove, ``Cost-optimal
    	scheduling in hybrid iaas clouds for deadline constrained workloads,'' in
    	\emph{Cloud Computing (CLOUD), 2010 IEEE 3rd International Conference
    		on}.\hskip 1em plus 0.5em minus 0.4em\relax IEEE, 2010, pp. 228--235.
    	
    	\bibitem{batista2007set}
    	D.~M. Batista, N.~L. Da~Fonseca, and F.~K. Miyazawa, ``A set of schedulers for
    	grid networks,'' in \emph{Proceedings of the 2007 ACM symposium on Applied
    		computing}.\hskip 1em plus 0.5em minus 0.4em\relax ACM, 2007, pp. 209--213.
    	
    	\bibitem{jiang2012joint}
    	J.~W. Jiang, T.~Lan, S.~Ha, M.~Chen, and M.~Chiang, ``Joint {VM} placement and
    	routing for data center traffic engineering,'' in \emph{INFOCOM}.\hskip 1em
    	plus 0.5em minus 0.4em\relax IEEE, 2012, pp. 2876--2880.
    	
    	\bibitem{calinescu2011maximizing}
    	G.~Calinescu, C.~Chekuri, M.~P{\'a}l, and J.~Vondr{\'a}k, ``Maximizing a
    	monotone submodular function subject to a matroid constraint,'' \emph{SIAM
    		Journal on Computing}, vol.~40, no.~6, pp. 1740--1766, 2011.
    	
    	\bibitem{filmus2012tight}
    	Y.~Filmus and J.~Ward, ``A tight combinatorial algorithm for submodular
    	maximization subject to a matroid constraint,'' in \emph{Foundations of
    		Computer Science (FOCS), 2012 IEEE 53rd Annual Symposium on}.\hskip 1em plus
    	0.5em minus 0.4em\relax IEEE, 2012, pp. 659--668.
    	
    	\bibitem{Bertsekas1989parallel}
    	D.~P. Bertsekas and J.~N. Tsitsiklis, \emph{Parallel and Distributed
    		Computation: Numerical Methods}.\hskip 1em plus 0.5em minus 0.4em\relax
    	Prentice Hall Englewood Cliffs, NJ, 1989, vol.~23.
    	
    	\bibitem{gill2016bidcache}
    	A.~S. Gill, L.~D'Acunto, K.~Trichias, and R.~van Brandenburg, ``Bid{C}ache:
    	{A}uction-based in-network caching in {ICN},'' in \emph{Globecom
    		Wkshps.}\hskip 1em plus 0.5em minus 0.4em\relax IEEE, 2016, pp. 1--6.
    	
    	\bibitem{xiao2004fast}
    	L.~Xiao and S.~Boyd, ``Fast linear iterations for distributed averaging,''
    	\emph{Systems \& Control Letters}, vol.~53, pp. 65--78, 2004.
    	
    	\bibitem{olfati2007consensus}
    	R.~Olfati-Saber, J.~A. Fax, and R.~M. Murray, ``Consensus and cooperation in
    	networked multi-agent systems,'' \emph{Proceedings of the IEEE}, vol.~95,
    	no.~1, pp. 215--233, 2007.
    	
    	\bibitem{gabber1981explicit}
    	O.~Gabber and Z.~Galil, ``Explicit constructions of linear-sized
    	superconcentrators,'' \emph{Journal of Computer and System Sciences},
    	vol.~22, no.~3, pp. 407--420, 1981.
    	
    	\bibitem{kleinberg2000small}
    	J.~Kleinberg, ``The small-world phenomenon: An algorithmic perspective,'' in
    	\emph{STOC}, 2000.
    	
    	\bibitem{watts1998collective}
    	D.~J. Watts and S.~H. Strogatz, ``Collective dynamics of `small-world'
    	networks,'' \emph{Nature}, vol. 393, no. 6684, pp. 440--442, 1998.
    	
    	\bibitem{barabasi1999emergence}
    	A.-L. Barab{\'a}si and R.~Albert, ``Emergence of scaling in random networks,''
    	\emph{Science}, vol. 286, no. 5439, pp. 509--512, 1999.
    	
    \end{thebibliography}

	
	\appendix
	\section{Appendix}
	
	\subsection{Consensus Algorithm} \label{secConsensusAlg}
	
	Here, we present an average consensus algorithm for all nodes to compute $\frac{1 }{|\mathcal{V}|}\sum_{v\in \mathcal{V}}C_{v0}$; see, e.g., \cite{xiao2004fast,olfati2007consensus}. Suppose each node $v\in \mathcal{V}$ initializes $s_v(0) = C_{v0}$ at time $t=0$ and updates $s_v$ according to a distributed linear iteration (involving only direct neighbor message exchanges):  
	\begin{equation*}
	s_v(t+1) = a_{vv}s_v(t)+ \sum_{u\in \mathcal{N}_v} a_{vu}s_u(t), \quad t=0,1,\ldots
	\end{equation*}
	where $[a_{uv}]$ satisfies one of the following conditions: 
	\begin{itemize}
		\item local-degree weights: 
		\begin{equation*}
		a_{vu} = \begin{cases}
		{1}/{\max\{ |\mathcal{N}_v|, |\mathcal{N}_u| \}} & (vu)\in \mathcal{E}\\
		1-\sum_{k\in \mathcal{N}_v} a_{vk} & u=v\\
		0 & \text{else}
		\end{cases}
		\end{equation*}
		
		\item constant edge weights: 
		\begin{equation*}
		a_{vu} = \begin{cases}
		\alpha & (vu)\in \mathcal{E}\\
		1- \alpha |\mathcal{N}_v| & u=v\\
		0 & \text{else}
		\end{cases}
		\end{equation*}
		with any $\alpha$ in the range $0<\alpha < \frac{2}{\max_{(vu)\in \mathcal{E}}  |\mathcal{N}_v| + |\mathcal{N}_u|  }$. 
	\end{itemize} 
	It follows that, if the network $\mathcal{G}$ is connected, then
	\begin{equation*}
	\lim_{t\to \infty} ~s_v(t) = \frac{\sum_{v\in \mathcal{V}}C_{v0}}{|\mathcal{V}|}, \quad \forall v\in \mathcal{V}.
	\end{equation*}
	Moreover, the convergence is \emph{geometric} at a  rate bounded above by the second largest eigenvalue of matrix $A = [a_{uv}]$.

\end{document}